\documentclass[10pt,headrule,fleqn,onecolumn]{article}

\usepackage{CJK} 
\usepackage{amssymb}
\usepackage{amsmath}
\usepackage{amsfonts}
\usepackage{amsopn}
\usepackage{amscd}
\usepackage{amsthm}
\usepackage{amsfonts,amssymb,color}
\usepackage{bm}
\usepackage{cite}
\usepackage{curves}
\usepackage{CJK}
\usepackage{delarray}
\usepackage{dsfont}
\usepackage{exscale}
\usepackage{framed}
\usepackage{hyperref}
\usepackage{indentfirst,latexsym,bm}
\usepackage{mathrsfs}
\usepackage{pifont}
\usepackage{xypic}

\usepackage{geometry}
\date{}
\begin{document}
\begin{CJK*}{GBK}{song}
\title{\bfseries\scshape Some existence theorems on all fractional $(g,f)$-factors with prescribed properties
}
\author{\small Sizhong Zhou$^{1}$\footnote{Corresponding
author. E-mail address: zsz\_cumt@163.com (S. Zhou)},  Tao Zhang$^{2}$\\
\small $1$. School of Mathematics and Physics, Jiangsu University of Science and Technology,\\
\small Mengxi Road 2, Zhenjiang, Jiangsu 212003, P. R. China\\
\small $2$. School of Economic and management, Jiangsu University of Science and Technology,\\
 \small Mengxi Road 2, Zhenjiang, Jiangsu 212003, P. R. China\\
}

\maketitle

\begin{abstract}
\noindent Let $G$ be a graph, and $g,f:V(G)\rightarrow Z^{+}$ with $g(x)\leq f(x)$ for each
$x\in V(G)$. We say that $G$ admits all fractional $(g,f)$-factors if $G$ contains an
fractional $r$-factor for every $r:V(G)\rightarrow Z^{+}$ with $g(x)\leq r(x)\leq f(x)$ for
any $x\in V(G)$. Let $H$ be a subgraph of $G$. We say that $G$ has all fractional
$(g,f)$-factors excluding $H$ if for every $r:V(G)\rightarrow Z^{+}$ with
$g(x)\leq r(x)\leq f(x)$ for all $x\in V(G)$, $G$ has a fractional $r$-factor $F_h$ such
that $E(H)\cap E(F_h)=\emptyset$, where $h:E(G)\rightarrow [0,1]$ is a function. In this
paper, we show a characterization for the existence of all fractional $(g,f)$-factors
excluding $H$ and obtain two sufficient conditions for a graph to have all fractional
$(g,f)$-factors excluding $H$.
\\
\begin{flushleft}
{\em Keywords:} graph; fractional $(g,f)$-factor; all fractional $(g,f)$-factors.

(2010) Mathematics Subject Classification: 05C70, 05C72
\end{flushleft}
\end{abstract}
\section{Introduction}
We consider finite undirected graphs without multiple edges or loops. Let $G=(V(G),E(G))$
be a graph, where $V(G)$ and $E(G)$ denote its vertex set and edge set, respectively. For any
$x\in V(G)$, $d_G(x)$ denotes the degree of $x$ in $G$. For any $S\subseteq V(G)$, the subgraph
of $G$ induced by $S$ is denoted by $G[S]$, and we write $G-S$ for $G[V(G)\setminus S]$. For two
disjoint vertex subsets $S$ and $T$ of $G$, we denote the number of edges with one end in $S$ and
the other end in $T$ by $e_G(S,T)$. The number of components of a graph $G$ is denoted by $\omega(G)$.

Let $g$ and $f$ be two positive integer-valued functions defined on $V(G)$ satisfying $g(x)\leq f(x)$
for any $x\in V(G)$. We define a $(g,f)$-factor as a spanning subgraph $F$ of $G$ which satisfies
$g(x)\leq d_F(x)\leq f(x)$ for any $x\in V(G)$. A $(g,f)$-factor is called an $f$-factor
if $g(x)=f(x)$ for all $x\in V(G)$. If $f(x)\equiv k$, then an $f$-factor is a $k$-factor, where $k$
is a positive integer. If $G$ admits an $r$-factor for every $r:V(G)\rightarrow Z^{+}$ which
satisfies $g(x)\leq r(x)\leq f(x)$ for any $x\in V(G)$ and $r(V(G))$ is even, then we say that $G$
has all $(g,f)$-factors.

Let $h:E(G)\rightarrow [0,1]$ be a function. For each $x\in V(G)$, $E(x)$ denotes the set of edges
incident with $x$. If $g(x)\leq\sum_{e\in E(x)}h(e)\leq f(x)$ holds for any $x\in V(G)$, then we
say that $F_h$ is an fractional $(g,f)$-factor of $G$ with indicator function $h$, where $F_h$ is
a subgraph with vertex set $V(G)$ and edge set $E_h=\{e:e\in E(G),h(e)>0\}$. An fractional
$(f,f)$-factor is said to be an fractional $f$-factor. An fractional $f$-factor is said to be an
fractional $k$-factor if $f(x)=k$ for each $x\in V(G)$. If $G$ includes an fractional $r$-factor
for every $r:V(G)\rightarrow Z^{+}$ with $g(x)\leq r(x)\leq f(x)$ for any $x\in V(G)$, then we say
that $G$ has all fractional $(g,f)$-factors. Let $H$ be a subgraph of $G$. If for every
$r:V(G)\rightarrow Z^{+}$ with $g(x)\leq r(x)\leq f(x)$ for any $x\in V(G)$, $G$ has an fractional
$r$-factor $F_h$ such that $E(H)\cap E(F_h)=\emptyset$, then we say that $G$ has all fractional
$(g,f)$-factors excluding $H$.

Let $\mathscr{F}$ be a set of graphs. If each component $C$ of $F$ is isomorphic to some member
of $\mathscr{F}$, then we say that $F$ is an $\mathscr{F}$-factor. A $\{K_n|n\geq1\}$-factor is
said to be a complete-factor.

Many authors have studied factors \cite{G,KO,DFGM,Z2,Z3}, fractional factors \cite{K,GLXZ,Z1,ZYS}
and all fractional factors \cite{Lu,ZS} in graphs. Lov\'asz \cite{L} present a characterization of
graphs having a $(g,f)$-factor.

\medskip

\noindent{\textbf{Theorem 1}} (Lov\'asz \cite{L}). Let $G$ be a graph, and let $g$ and $f$ be two
integer-valued functions defined on $V(G)$ with $0\leq g(x)\leq f(x)$ for any $x\in V(G)$. Then $G$
admits a $(g,f)$-factor if and only if
$$
f(S)+\sum_{x\in T}d_{G-S}(x)-g(T)-h_G(S,T,g,f)\geq0
$$
for all disjoint subsets $S$ and $T$ of $V(G)$, where $f(S)=\sum_{x\in S}f(x)$, $g(T)=\sum_{x\in T}g(x)$
and $h_G(S,T,g,f)$ denotes the number of components
$C$ of $G-(S\cup T)$ with $g(x)=f(x)$ for any $x\in V(C)$ and $f(V(C))+e_G(V(C),T)\equiv1 \ (mod \ 2)$.

\medskip

Niessen \cite{N} showed a result for a graph to have all $(g,f)$-factors.

\medskip

\noindent{\textbf{Theorem 2}} (Niessen \cite{N}). A graph $G$ has all $(g,f)$-factors if and only if
$$
g(S)+\sum_{x\in T}d_{G-S}(x)-f(T)-q_G(S,T,g,f)\geq\left\{
\begin{array}{ll}
-1,&if \ f\neq g\\
0,&if \ f=g\\
\end{array}
\right.
$$
for all disjoint subsets $S$ and $T$ of $V(G)$, where $q_G(S,T,g,f)$
denotes the number of components $C$ of $G-(S\cup T)$ such that
there exists a vertex $v\in V(C)$ with $g(v)<f(v)$ or
$e_G(V(C),T)+f(V(C))\equiv1$ $(mod \ 2)$.

\medskip

Anstee \cite{A} posed a necessary and sufficient condition for a graph to have a fractional
$(g,f)$-factor. Liu and Zhang \cite{LZ} gave a new proof.

\medskip

\noindent{\textbf{Theorem 3}} (Anstee \cite{A}, Liu and Zhang \cite{LZ}). Let $G$ be a graph,
and $g,f:V(G)\rightarrow Z^{+}$ be two integer-valued functions with $g(x)\leq f(x)$ for any
$x\in V(G)$. Then $G$ admits a fractional $(g,f)$-factor if and only if
$$
f(S)+\sum_{x\in T}d_{G-S}(x)-g(T)\geq0
$$
for any subset $S$ of $V(G)$, where $T=\{x:x\in V(G)-S, d_{G-S}(x)<g(x)\}$.

\medskip

Lu \cite{Lu} obtained a result similar to Theorem 2 for all fractional $(g,f)$-factors.

\medskip

\noindent{\textbf{Theorem 4}} (Lu \cite{Lu}). Let $G$ be a graph and $g,f:V(G)\rightarrow Z^{+}$
be two integer-valued functions with $g(x)\leq f(x)$ for any $x\in V(G)$. Then $G$ has all
fractional $(g,f)$-factors if and only if
$$
g(S)+\sum_{x\in T}d_{G-S}(x)-f(T)\geq0
$$
for any subset $S$ of $V(G)$, where $T=\{x:x\in V(G)-S, d_{G-S}(x)<f(x)\}$.

\medskip

In this paper, we investigate the existence of all fractional $(g,f)$-factors excluding any
given subgraph in graphs, and obtain three new results which are shown in the following
section.

\section{Main results}
In this section, we show our main results.

\medskip

\noindent{\textbf{Theorem 5}}. Let $G$ be a graph and $g,f:V(G)\rightarrow Z^{+}$ be two
integer-valued functions with $g(x)\leq f(x)$ for each $x\in V(G)$. Let $H$ be a subgraph of
$G$. Then $G$ admits all fractional $(g,f)$-factors excluding $H$ if and only if
$$
g(S)+\sum_{x\in T}d_{G-S}(x)-f(T)\geq\sum_{x\in T}d_H(x)-e_H(S,T)
$$
for any subset $S$ of $V(G)$, where $T=\{x:x\in V(G)-S, d_{G-S}(x)-d_H(x)+e_H(x,S)<f(x)\}$.

\medskip

If $E(H)=\emptyset$ in Theorem 5, then Theorem 4 is obtained immediately. Therefore, Theorem 4
is a special case of Theorem 5. By using Theorem 5, we obtain two sufficient conditions for
a graph to have all fractional $(g,f)$-factors excluding $H$, which are the following theorems.

\medskip

\noindent{\textbf{Theorem 6}}. Let $G$ be a graph, $F$ be a complete-factor of $G$
with $\omega(F)\geq2$ and $g,f:V(G)\rightarrow Z^{+}$ be two integer-valued functions
with $g(x)\leq f(x)$ for any $x\in V(G)$. Let $H$ be a subgraph of $G$. If $G-V(C)$ has all
fractional $(g,f)$-factors excluding $H$ for each component $C$ of $F$, then $G$
itself admits all fractional $(g,f)$-factors excluding $H$.

\medskip

If $E(H)=\emptyset$ in Theorem 6, then we obtain the following corollary.

\medskip

\noindent{\textbf{Corollary 7}}. Let $G$ be a graph, $F$ be a complete-factor of $G$
with $\omega(F)\geq2$ and $g,f:V(G)\rightarrow Z^{+}$ be two integer-valued functions
with $g(x)\leq f(x)$ for any $x\in V(G)$. If $G-V(C)$ has all fractional $(g,f)$-factors for
each component $C$ of $F$, then $G$ itself admits all fractional $(g,f)$-factors.

\medskip

If $F$ is a 1-factor of $G$ in Theorem 6, then we have the following corollary.

\medskip

\noindent{\textbf{Corollary 8}}. Let $G$ be a graph, $F$ be a 1-factor of $G$ and
$g,f:V(G)\rightarrow Z^{+}$ be two integer-valued functions with $g(x)\leq f(x)$ for any
$x\in V(G)$. Let $H$ be a subgraph of $G$. If $G-\{x,y\}$ has all fractional $(g,f)$-factors
excluding $H$ for each $xy\in E(F)$, then $G$ itself has all fractional $(g,f)$-factors
excluding $H$.

\medskip

\noindent{\textbf{Theorem 9}}. Let $G$ be a graph, $g,f:V(G)\rightarrow Z^{+}$ be two
integer-valued functions with $g(x)\leq f(x)$ for any $x\in V(G)$, and $H$ be a subgraph of
$G$. If $d_G(x)\geq f(x)+d_H(x)$ and $g(x)(d_G(y)-d_H(y))\geq d_G(x)f(y)$ for any $x,y\in V(G)$,
then $G$ has all fractional $(g,f)$-factors excluding $H$.

\section{The proofs of Main Theorems}
\noindent{\it Proof of Theorem 5.} \ We first verify sufficiency. Let $r:V(G)\rightarrow Z^{+}$
be an arbitrary integer-valued function with $g(x)\leq r(x)\leq f(x)$ for any $x\in V(G)$. By
the definition of all fractional $(g,f)$-factors excluding $H$, we need only to prove that $G$
has a fractional $r$-factor excluding $H$, that is, we need only to verify that $G'$ has a
fractional $r$-factor, where $G'=G-E(H)$.

For any $S\subseteq V(G)$, we write $T=\{x:x\in V(G)-S,d_{G'-S}(x)<f(x)\}$ and
$T'=\{x:x\in V(G')-S,d_{G'-S}(x)<r(x)\}$. Note that
$T=\{x:x\in V(G)-S,d_{G-S}(x)-d_H(x)+e_H(x,S)<f(x)\}$. Thus, we obtain
\begin{eqnarray*}
&&r(S)+\sum_{x\in T'}d_{G'-S}(x)-r(T')\geq g(S)+\sum_{x\in T}d_{G'-S}(x)-f(T)\\
&&=g(S)+\sum_{x\in T}d_{G-S}(x)-f(T)-\sum_{x\in T}d_H(x)+e_H(S,T)\geq0.
\end{eqnarray*}
According to Theorem 3, $G'$ admits a fractional $r$-factor.

Now we verify the necessary. Conversely, we assume that there exists some subset $S$ of $V(G)$
satisfying
$$
g(S)+\sum_{x\in T}d_{G-S}(x)-f(T)<\sum_{x\in T}d_H(x)-e_H(S,T),
$$
where $T=\{x:x\in V(G)-S,d_{G-S}(x)-d_H(x)+e_H(x,S)<f(x)\}$. Set $G'=G-E(H)$,
$r(x)=g(x)$ for any $x\in S$ and $r(y)=f(y)$ for any $y\in V(G)\setminus S$. Obviously,
$T=\{x:x\in V(G')-S,d_{G'-S}(x)<r(x)$. Thus, we obtain
\begin{eqnarray*}
0&>&g(S)+\sum_{x\in T}d_{G-S}(x)-f(T)-\sum_{x\in T}d_H(x)+e_H(S,T)\\
&=&r(S)+\sum_{x\in T}d_{G'-S}(x)-r(T).
\end{eqnarray*}
In terms of Theorem 3, $G'$ has no fractional $r$-factor. And so, $G'$ has no all fractional
$(g,f)$-factors. Combining this with $G'=G-E(H)$, $G$ has no all fractional $(g,f)$-factors
excluding $H$, a contradiction. Theorem 5 is proved. \hfill $\Box$

\medskip

\medskip

\noindent{\it Proof of Theorem 6.} \ For any $S\subseteq V(G)$, we write
$$
\delta_G(S,T)=g(S)+\sum_{x\in T}d_{G-S}(x)-f(T),
$$
where $T=\{x:x\in V(G)-S,d_{G-S}(x)-d_H(x)+e_H(x,S)<f(x)\}$. According to Theorem 5, in order to
prove the theorem we need only to show that
$$
\delta_G(S,T)\geq\sum_{x\in T}d_H(x)-e_H(S,T)
$$
for any $S\subseteq V(G)$ and $T=\{x:x\in V(G)-S,d_{G-S}(x)-d_H(x)+e_H(x,S)<f(x)\}$.

Set $U=V(G)-(S\cup T)$. For each component $C$ of $F$, we write $S_C=V(C)\cap S$, $T_C=V(C)\cap T$
and $U_C=V(C)\cap U$. It is easy to see that
\begin{eqnarray*}
&&\sum_{x\in(T-V(C))}d_{(G-V(C))-(S-V(C))}(x)\\
&=&\sum_{x\in(T-V(C))}d_{G-V(C)}(x)-e_{G-V(C)}(T-V(C),S-V(C))\\
&=&\sum_{x\in(T-V(C))}d_G(x)-e_G(V(C),T-V(C))-e_G(T-V(C),S-V(C))\\
&=&\sum_{x\in(T-V(C))}d_{G-S}(x)+e_G(T-V(C),S)-e_G(T-V(C),V(C))\\
&&-e_G(T-V(C),S-V(C))\\
&=&\sum_{x\in(T-V(C))}d_{G-S}(x)+e_G(T-V(C),S)-e_G(T-V(C),T_C\cup S_C\cup U_C)\\
&&-e_G(T-V(C),S-V(C))\\
&=&\sum_{x\in(T-V(C))}d_{G-S}(x)+e_G(T-V(C),S)-e_G(T-V(C),T_C)-e_G(T-V(C),S_C)\\
&&-e_G(T-V(C),U_C)-e_G(T-V(C),S)+e_G(T-V(C),S_C)\\
&=&\sum_{x\in(T-V(C))}d_{G-S}(x)-e_G(T-V(C),T_C)-e_G(T-V(C),U_C)
\end{eqnarray*}
and
\begin{eqnarray*}
&&\sum_{x\in(T-V(C))}d_H(x)-e_H(S-V(C),T-V(C))\\
&=&\sum_{x\in T}d_H(x)-\sum_{x\in T_C}d_H(x)-e_H(S-V(C),T)+e_H(S-V(C),T_C)\\
&=&\sum_{x\in T}d_H(x)-\sum_{x\in T_C}d_H(x)-e_H(S,T)+e_H(S_C,T)+e_H(S,T_C)-e_H(S_C,T_C)\\
&\geq&\sum_{x\in T}d_H(x)-e_H(S,T)-\sum_{x\in T_C}d_H(x)+e_H(S,T_C)
\end{eqnarray*}
for each component $C$ of $F$. Since $G-V(C)$ has all fractional $(g,f)$-factors excluding $H$
for each component $C$ of $F$, it follows from Theorem 5 that
\begin{eqnarray*}
&&\sum_{x\in T}d_H(x)-e_H(S,T)-\sum_{x\in T_C}d_H(x)+e_H(S,T_C)\\
&\leq&\sum_{x\in(T-V(C))}d_H(x)-e_H(S-V(C),T-V(C))\\
&\leq&\delta_{G-V(C)}(S-V(C),T-V(C))\\
&=&g(S-V(C))+\sum_{x\in(T-V(C))}d_{(G-V(C)-(S-V(C))}(x)-f(T-V(C))\\
&=&g(S-V(C))+\sum_{x\in(T-V(C))}d_{G-S}(x)-e_G(T-V(C),T_C)-e_G(T-V(C),U_C)-f(T-V(C))\\
&\leq&g(S-V(C))+\sum_{x\in(T-V(C))}d_{G-S}(x)-f(T-V(C))\\
&=&g(S)-g(S_C)+\sum_{x\in T}d_{G-S}(x)-\sum_{x\in T_C}d_{G-S}(x)-f(T)+f(T_C)\\
&=&\delta_G(S,T)-g(S_C)-\sum_{x\in T_C}d_{G-S}(x)+f(T_C),
\end{eqnarray*}
which implies
$$
\delta_G(S,T)-\Big(\sum_{x\in T}d_H(x)-e_H(S,T)\Big)\geq g(S_C)+\sum_{x\in T_C}d_{G-S}(x)-f(T_C)-\Big(\sum_{x\in T_C}d_H(x)-e_H(S,T_C)\Big),
$$
and so
\begin{eqnarray*}
&&\sum_{j=1}^{\omega(F)}\delta_G(S,T)-\sum_{j=1}^{\omega(F)}\Big(\sum_{x\in T}d_H(x)-e_H(S,T)\Big)\\
&\geq&\sum_{j=1}^{\omega(F)}\Big(g(S_C)+\sum_{x\in T_C}d_{G-S}(x)-f(T_C)\Big)-\sum_{j=1}^{\omega(F)}\Big(\sum_{x\in T_C}d_H(x)-e_H(S,T_C)\Big)\\
&=&g(S)+\sum_{x\in T}d_{G-S}(x)-f(T)-\Big(\sum_{x\in T}d_H(x)-e_H(S,T)\Big)\\
&=&\delta_G(S,T)-\Big(\sum_{x\in T}d_H(x)-e_H(S,T)\Big).
\end{eqnarray*}
Thus, we obtain
$$
(\omega(F)-1)\delta_G(S,T)\geq(\omega(F)-1)\Big(\sum_{x\in T}d_H(x)-e_H(S,T)\Big).
$$
Note that $\omega(F)\geq2$. Hence, we have
$$
\delta_G(S,T)\geq\sum_{x\in T}d_H(x)-e_H(S,T).
$$
By Theorem 5, $G$ has all fractional $(g,f)$-factors excluding $H$. This completes
the proof of Theorem 6. \hfill $\Box$

\medskip

\medskip

\noindent{\it Proof of Theorem 9.} \ In terms of Theorem 5, we need only to prove that
$$
g(S)+\sum_{x\in T}d_{G-S}(x)-f(T)\geq\sum_{x\in T}d_H(x)-e_H(S,T)
$$
for any subset $S$ of $V(G)$, where $T=\{x:x\in V(G)-S, d_{G-S}(x)-d_H(x)+e_H(x,S)<f(x)\}$.

If $S=\emptyset$, then by $d_G(x)\geq f(x)+d_H(x)$ for any $x\in V(G)$,
$$
g(S)+\sum_{x\in T}d_{G-S}(x)-f(T)=\sum_{x\in T}d_G(x)-f(T)\geq\sum_{x\in T}d_H(x)-e_H(S,T).
$$

In the following, we consider $S\neq\emptyset$. Note that $g(x)(d_G(y)-d_H(y))\geq d_G(x)f(y)$
for any $x,y\in V(G)$. Thus, we obtain
$$
\Big(\sum_{x\in S}d_G(x)\Big)f(y)\leq g(S)\Big(d_G(y)-d_H(y)\Big)
$$
and
$$
\Big(\sum_{x\in S}d_G(x)\Big)f(T)\leq g(S)\Big(\sum_{x\in T}d_G(y)-\sum_{x\in T}d_H(y)\Big).\eqno(1)
$$
Note that
$$
\sum_{x\in S}d_G(x)-\sum_{x\in T}d_G(y)\geq-\sum_{x\in T}d_{G-S}(x).\eqno(2)
$$
It follows from (1) and (2) that
\begin{eqnarray*}
&&\Big(\sum_{x\in S}d_G(x)\Big)\Big(g(S)+\sum_{x\in T}d_{G-S}(x)-f(T)-\sum_{x\in T}d_H(x)+e_H(S,T)\Big)\\
&&=\Big(\sum_{x\in S}d_G(x)\Big)g(S)+\Big(\sum_{x\in S}d_G(x)\Big)\Big(\sum_{x\in T}d_{G-S}(x)\Big)-\Big(\sum_{x\in S}d_G(x)\Big)f(T)\\
&& \ \ \ -\Big(\sum_{x\in S}d_G(x)\Big)\Big(\sum_{x\in T}d_H(x)\Big)+\Big(\sum_{x\in S}d_G(x)\Big)e_H(S,T)\\
&&\geq\Big(\sum_{x\in S}d_G(x)\Big)g(S)+\Big(\sum_{x\in S}d_G(x)\Big)\Big(\sum_{x\in T}d_{G-S}(x)\Big)-g(S)\Big(\sum_{x\in T}d_G(y)-\sum_{x\in T}d_H(y)\Big)\\
&& \ \ \ -\Big(\sum_{x\in S}d_G(x)\Big)\Big(\sum_{x\in T}d_H(x)\Big)+\Big(\sum_{x\in S}d_G(x)\Big)e_H(S,T)\\
&&=g(S)\Big(\sum_{x\in S}d_G(x)-\sum_{x\in T}d_G(y)\Big)+\Big(\sum_{x\in S}d_G(x)\Big)\Big(\sum_{x\in T}d_{G-S}(x)\Big)+g(S)\Big(\sum_{x\in T}d_H(y)\Big)\\
&& \ \ \ -\Big(\sum_{x\in S}d_G(x)\Big)\Big(\sum_{x\in T}d_H(x)\Big)+\Big(\sum_{x\in S}d_G(x)\Big)e_H(S,T)\\
&&\geq-g(S)\Big(\sum_{x\in T}d_{G-S}(x)\Big)+\Big(\sum_{x\in S}d_G(x)\Big)\Big(\sum_{x\in T}d_{G-S}(x)\Big)+g(S)\Big(\sum_{x\in T}d_H(y)\Big)\\
&& \ \ \ -\Big(\sum_{x\in S}d_G(x)\Big)\Big(\sum_{x\in T}d_H(x)\Big)+\Big(\sum_{x\in S}d_G(x)\Big)e_H(S,T)\\
&&=\Big(\sum_{x\in S}d_G(x)-g(S)\Big)\Big(\sum_{x\in T}d_{G-S}(x)-\sum_{x\in T}d_H(y)\Big)+\Big(\sum_{x\in S}d_G(x)\Big)e_H(S,T)\\
&&\geq\Big(\sum_{x\in S}d_G(x)-g(S)\Big)(-e_H(S,T))+\Big(\sum_{x\in S}d_G(x)\Big)e_H(S,T)\\
&&=g(S)e_H(S,T)\geq0.
\end{eqnarray*}
It is obvious that $\sum_{x\in S}d_G(x)\geq f(S)\geq g(S)\geq|S|\geq1$. Hence, we have
$$
g(S)+\sum_{x\in T}d_{G-S}(x)-f(T)-\sum_{x\in T}d_H(x)+e_H(S,T)\geq0,
$$
that is,
$$
g(S)+\sum_{x\in T}d_{G-S}(x)-f(T)\geq\sum_{x\in T}d_H(x)-e_H(S,T).
$$
This finishes the proof of Theorem 9. \hfill $\Box$

\medskip

{\bf Acknowledgments.} The authors would like to thank the anonymous
referees for their kind help and valuable suggestions which led to
an improvement of this paper.

\end{CJK*}
\end{document}